\documentclass[12pt,a4paper,reqno]{amsart}

\usepackage[british]{babel}
\usepackage{amsmath,amsthm,amssymb,mathtools}
\usepackage[alphabetic,abbrev, msc-links]{amsrefs}
\usepackage{enumitem,lineno}
\usepackage[normalem]{ulem}
\usepackage{hyperref}
\usepackage{esvect}
\usepackage[nameinlink]{cleveref}
\usepackage{geometry}

\usepackage{xcolor}

\newcommand\numberthis{\addtocounter{equation}{1}\tag{\theequation}}
\newcommand{\braces}[1]{\left\lbrace{#1}\right\rbrace}
\newcommand{\verts}[1]{\left\lvert {#1} \right\rvert}
\newcommand{\mbbN}{\mathbb{N}}
\newcommand{\Kalephzero}{K^{\aleph_0}}
\newcommand{\gm}{G_M}
\newcommand{\nm}{N_M}
\newcommand{\gb}{G_B}

\newcommand{\thiscolour}{c_p}

\newtheorem{theorem}{Theorem}\numberwithin{theorem}{section}
\newtheorem{lem}[theorem]{Lemma}

\theoremstyle{definition}
\newtheorem{question}[theorem]{Question}

\DeclareMathOperator{\im}{im}

\title[Infinite Maker-Breaker-Games]{The $K^{\aleph_0}$ game: Vertex colouring}
\author[Bowler, Emde \& Gut]{Nathan Bowler, Marit Emde \& Florian Gut}
\address{Universität Hamburg, Faculty of Mathematics, Informatics and Natural Sciences, Department of Mathematics, Bundesstra{\ss}e 55 (Geomatikum), 20146 Hamburg, Germany}
\email{nathan.bowler@uni-hamburg.de, florian.gut@uni-hamburg.de}

\begin{document}
	
	\begin{abstract}
		We investigate games played between Maker and Breaker on an infinite complete graph whose vertices are coloured with colours from a given set, each colour appearing infinitely often.
		The players alternately claim edges, Maker’s aim being to claim all edges of a sufficiently colourful infinite complete subgraph and Breaker’s aim being to prevent this.
		We show that if there are only finitely many colours then Maker can obtain a complete subgraph in which all colours appear infinitely often, but that Breaker can prevent this if there are infinitely many colours.
		Even when there are infinitely many colours, we show that Maker can obtain a complete subgraph in which infinitely many of the colours each appear infinitely often.
	\end{abstract}
	
	\maketitle
	
	\section{Introduction}\label{sec:Introduction}
	Games have been of interest to mathematicians for centuries.
	The field was sparked by the analysis of historic games such as tic-tac-toe.
	This is a finite game and naturally, this was the first class of games analysed.
	In recent decades infinite games have increasingly drawn the attention of researchers.
	An intuitive  starting point is just moving the rules of a finite game to an infinite board.
	Consider the aforementioned Tic-tac-toe.
	Its counterpart on an infinite board became known as \emph{unrestricted $3$-in-a-row} and was further generalised to $n$\emph{-in-a-row} (Beck~\cite{B08}).
	As some results in the finite version simply rely on the fact that there are only finitely many possible plays, this often yields interesting insights.
	These types of games are called \emph{semi-infinite}; a comprehensive overview of the known results about such games can be found in~\cite{B08}.
	
	Games also give rise to an interesting field in set theory, which is e.g. described in (Moschovakis~\cite{M09}*{Chapter 6}):
	One may assume as an axiom that every game is determined, i.e. that for any two player game with complete information at least one of the players has a winning strategy.
	It is known that this is not consistent with the axiom of choice, but it nevertheless provides a good framework to study descriptive set theory, itself an essential field of research in present-day mathematical logic.
	
	However, there has not yet been a systematic analysis of infinite combinatorial games, i.e. games played on an infinite board with perfect information where the players do moves sequentially.
	One classic variant is the following two player game, which we will call the {\em strong  $H$-building game}.
	The game depends on a fixed graph $H$.
	During the course of the game the players alternately claim edges of a sufficiently large complete graph $K^n$.
	The game ends as soon as $H$ is contained as a subgraph in the graph induced by some player's claimed edges.
	The word {\em strong} in the name of the game indicates that  both players have the same aim (to build a copy of $H$), thus the two players' roles differ only in who plays first. This is to distinguish it from the {\em Maker-Breaker} variant of the game, in which only one of the players (Maker) is trying to build a copy of $H$ and the other (Breaker) is simply trying to prevent this. In particular, Breaker does not win simply by virtue of building his own copy of $H$ first.
	
	Finite strong games have been extensively studied, and in particular the strong $H$-building game is fairly well understood.  Indeed, as long as the board is sufficiently large, any such game is a first player win:
	due to Ramsey's Theorem (see Diestel~\cite{D17}*{Theorem 9.1.1}) after all edges of the $K^n$ have been claimed, there will always be a copy of $H$ contained in one of the players' graphs if $n \geq R(\vert H \vert )$.
	So one of the players must have a winning strategy, and using a technique called \emph{strategy stealing} (see Hefetz, Krivelevich, Stojakovi\'c and Szab\'o~\cite{HKSS14}*{Theorem 1.3.1}) it cannot be the second player.
	The argument can be roughly sketched as follows:
	suppose the second player has a winning strategy.
	Then the first player can make an arbitrary move in his first turn and from then on play according to the winning strategy as if she is the second player.
	
	Note that even though we now know that there is a winning strategy for the first player if $n \geq R(\vert H \vert )$, we only have an abstract argument and can not deduce an actual strategy according to which the first player should play.
	In particular, this argument does not show the existence of an upper bound only depending on $H$ for the number of moves the first player needs to win.
	
	The problem of establishing such bounds is in general very hard, and for example if $H$ is a complete graph this problem has only been resolved when $H$ has at most 4 vertices (see Beck~\cite{B02}).
	By contrast, such bounds are well known for the Maker-Breaker $H$-building game (see~\cite{HKSS14}*{Chapter 2}).
	Indeed, Maker-Breaker variants are usually easier to analyse, and are often investigated as a preliminary step before analysing the strong game. 
	
	In contrast to the fertility of the study of finite combinatorial games, infinite combinatorial games have thus far proved barren ground, in that the strong games are too hard to analyse and the Maker-Breaker games are too easy.
	Let us consider first the strong $H$-building game, but now played on an infinite complete graph.
	Since there are now infinitely many edges on the board to choose from these edges need not be exhausted in the course of the game, and there is no guarantee that the players will ever even claim all edges of a finite complete subgraph $K^n$ between them.
	Thus Ramsey's Theorem is no longer applicable.
	Play could continue forever without either player ever winning!
	
	So even though we can still use a strategy-stealing argument to rule out the existence of a winning strategy for the second player, it is possible that he has a strategy to force a draw. Although it might seem implausible that this could really happen, in fact for the analogous games played on 5-regular hypergraphs rather than graphs cases are known in which the second player can force a draw.
	Such an example was constructed by Hefetz, Kusch, Narins, Pokrovskiy, Requil\'{e} and Sarid in~\cite{HKNPRS17}.
	
	In fact, it is folklore that the existence of a winning strategy for the first player in the infinite strong $H$-building game is equivalent to the existence of a finite upper bound on the number of moves the first player needs to win in the finite strong $H$-building games, a straightforward proof by a compactness argument can be found in (Leader~\cite{L08}*{Proposition 4}), and we saw above that such problems currently seem intractable. 
	
	Then again, the Maker-Breaker $H$-building game is trivial; the winning strategies for the finite variants also work in the infinite variant.
	However, Erde recently noticed that there are interesting Maker-Breaker $H$-building games on infinite graphs; one just has to take $H$ itself to be infinite~\cite{E19}. 
	
	We believe that the investigation of such infinite Maker-Breaker games will prove very fruitful. In this paper we begin that investigation by considering a few simple variants on the $\Kalephzero$-game.
	We will present a winning strategy for Maker in the basic version of the game in \Cref{section_TheOriginalKOmegaGame}.
	
	A natural variation arises if we colour the vertices of the board beforehand and demand that Maker respects this colouring in such a way that her $\Kalephzero$ again contains infinitely many vertices of every colour.
	
	In \Cref{section:finitely_many_colours} we will present a winning strategy for Maker if there are only finitely many colours.
	If there are infinitely many colours then, on the one hand, Maker can still incorporate infinitely many different colours infinitely often into their $\Kalephzero$, which we will prove in \Cref{subsec:infinitely_many_vertex_colourclasses_Maker}.
	On the other hand, in case Makers' aim is to incorporate all colours that are present on the board, Breaker can stop Maker from doing so, which we will prove in \Cref{subsec:infinitely_many_vertex_colourclasses_Breaker_strategy} by giving a winning strategy for Breaker.
	In fact, we prove even more: Breaker can stop Maker from incorporating cofinitely many colours infinitely often into her $\Kalephzero$.
	
	Another possible variant of the game is to demand that Maker respect some predefined colouring of the edges of the board instead, we will pose questions about this and other problems in \Cref{sec:open_problems}.
	
	There are two new kinds of difficulty which arise for infinite games of this kind but which are not so relevant for finite games. First of all, infinite games need not have a winning strategy for either player.
	Indeed, assuming the continuum hypothesis there are undetermined infinite combinatorial Maker-Breaker games on an infinite board.
	We present such a game in \Cref{sec:Appendix}.
	
	Another difficulty arising for the games studied in this paper is that it would be hopeless for Maker to attempt to build her copy of $\Kalephzero$ up iteratively, by building larger and larger nested $K^i$.
	The difficulty is that for any fixed $K^i$ with $i\geq 2$ it is straightforward for Breaker to claim at least one edge to it from each subsequent vertex, thus preventing it from being integrated into any $\Kalephzero$ claimed by Maker.
	Thus Maker must proceed more speculatively, never knowing during the course of the game which vertices will be integrated into her eventual $\Kalephzero$.
	
	Let us now begin by introducing some necessary notation.
	
	\section{Preliminaries}\label{sec:Preliminaries}
	Throughout the paper we draw on the standard definitions as established in \cite{D17}, if not explicitly mentioned otherwise.
	There are two players, who alternately take \emph{turns} throughout the game.
	In each individual turn they pick an edge from the \emph{board} $G = \left(V,E\right) := \left(  \mbbN , \mbbN^{2} \right) \cong \Kalephzero$, and colour it in their respective colour.
	We will refer to Maker by \emph{she} or \emph{her} and likewise we refer to Breaker by \emph{he} or \emph{him}.
	We will assume that Maker colours her edges in magenta, which we will abbreviate with $M$ and that Breaker uses the colour blue which we will abbreviate with $B$.
	Neither player may colour an edge that was already picked by either player in a previous turn.
	The goal for Maker will be to colour the edges in such a way that her subgraph contains a specific substructure, i.e. a $\Kalephzero$.
	Breaker's goal is to stop Maker from doing so.
	
	For a colour $\gamma \in \braces{M,B}$ and for any point in the game we define $E(G_{\gamma})$ to be the edges that have been coloured in the colour $\gamma$ up to that point, $V(G_{\gamma})$ to be all the vertices that are incident with at least one edge of $E(G_{\gamma})$ and thus define the graph $G_{\gamma} = (V(G_{\gamma}), E(G_{\gamma}))$.
	$N_{\gamma}(v)$ are the \emph{neighbours} of a vertex $v \in V(G_{\gamma})$ in $G_{\gamma}$ and the $\gamma$-\emph{degree} $\deg_{\gamma}(v)$ of a vertex $v \in G$ is $\deg_{\gamma}(v) = \left\vert N_{\gamma}(v) \right\vert$.
	Accordingly, we will say that two vertices $v$ and $w$ are $\gamma$-connected or $\gamma$-adjacent, if $vw \in G_{\gamma}$.
	
	When we say that a player $\gamma \in \braces{M,B}$ $\gamma$-\emph{connects} a vertex $v$ to a vertex $w$ in a turn, then we mean that Maker or Breaker claims the edge $vw$ in that turn respectively.
	We will mean the same when we say that Maker or Breaker \emph{plays from} a vertex $v$ \emph{to} a vertex $w$.
	
	When we talk about a \emph{fresh} vertex, we mean a vertex $v \in V \setminus \left(  V(\gm) \cup V(\gb) \right)$.
	Since Maker adds vertices to $G_M$ in distinct moves, the vertices become ordered in a natural way.
	We will make use of that and assign indices to the vertices accordingly, i.e. $v_k$ is the $k$\textsuperscript{th} vertex that Maker adds to her subgraph $\gm$.
	When Maker claims an edge incident with two fresh vertices, she assigns the next two indices to these vertices arbitrarily, which will only happen on Maker's first turn in our construction.
	For natural numbers $n \in \mbbN $ we set $[n] := \braces{1,2,3,\dots,n}$.
	When we want to prove that a game is a win for Breaker, we shall always do this by means of a \emph{pairing strategy}. That is, we will define a family of disjoint pairs of edges from $E(G)$ with the intention that whenever Maker claims one edge from a pair Breaker claims the other one in his following turn.
	
	It will then suffice to verify, for the game in question, that any $\Kalephzero$ of the kind that Maker is trying to build must include both edges of at least one such pair.
	
	\section{The basic version}\label{section_TheOriginalKOmegaGame}
	We will begin by investigating the basic version of the $\Kalephzero$-game.
	In this game the aim of Maker is that $\gm$ contains a $\Kalephzero$ at the end of the game.
	We will prove that Maker can win this game.
	We will achieve this by first describing a strategy according to which Maker should play and then verifying that, in fact, $\Kalephzero \subseteq \gm$ holds true.
	
	Our focus will be on two different kinds of activity by Maker. 
	On the one hand, she will regularly want to add fresh vertices to her subgraph $\gm$.
	On the other hand, she must ensure that $\gm$ is as interconnected as possible and thus contains large complete graphs.
	The same interplay between making $\gm$ highly interconnected and regularly moving on to fresh vertices will also provide the basic rhythm for our strategies for Maker in later sections. 
	
	\medskip
	We will call the following strategy for Maker the \emph{structured greedy strategy}.
	In her first turn, she picks some edge $v_1v_2$.
	In case Breaker was the first player, she picks one that only uses fresh vertices.\\
	In a later turn, suppose $v_n$ is the last vertex that was added to Maker's subgraph.
	Now, if there is some $v_i$, $1 \leq i < n$ such that
	\begin{enumerate}[label={($\square${\arabic{*}})}]
		\item\label{vanilla_version_not_yet_claimed} $v_iv_n$ has not yet been claimed in either colour,
		\item\label{vanilla_version_neighbourhoods_match} $N_M (v_n) \subseteq N_M (v_i) $, and
		\item\label{vanilla_version_smallest_vertex} $i$ is minimal subject to~\labelcref{vanilla_version_not_yet_claimed,vanilla_version_neighbourhoods_match},
	\end{enumerate}
	then Maker claims $v_iv_n$.
	If there is no such $v_i$, she picks a fresh vertex $v_{n+1}$ and claims $v_1v_{n+1}$.
	
	\begin{theorem}\label{theo:Original_Game}
		The structured greedy strategy is a winning strategy for Maker in the basic version of the $\Kalephzero$-game.
	\end{theorem}
	\begin{proof}
		We consider an arbitrary play of the game in which Maker follows the structured greedy strategy.
		We must show that at the end of the game $\gm$ includes a $\Kalephzero$.
		We will recursively construct a complete graph $K^n \subseteq G_M$ as well as a set of vertices $W_n \subseteq V(G_M)$ for every $n \in \mbbN\setminus \braces{0}$ such that:
		\begin{enumerate}[label={($\blacksquare${\arabic{*}})}]
			\item\label{item:vanilla_version_verification_size} $\vert K^n \vert = n$, $\vert W_n \vert = \aleph_0$,
			\item\label{item:vanilla_version_verification_inclusion} $K^n \subseteq K^{n+1}$ for $n > 1$, and
			\item\label{item:vanilla_version_verification_connected_in_time} for any $w \in W_n$ the first $n$ vertices to which $w$ was $M$-connected were $V(K^n)$.
		\end{enumerate}
		If we successfully construct such a sequence $K^1 \subset K^2 \subset K^3 \subset \dots$, the claim follows immediately for
		\[
		\bigcup_{i \in \mbbN} K^i = \Kalephzero .
		\]
		The purpose of the sets $W_n$ is to ensure that there will be a suitable candidate to enlarge the complete graph at each step.
		
		\textbf{Initial step}:
		We can set $K^1 = \left( \braces{v_1}, \emptyset \right) $ and $W_1 = V(G_M) \setminus \braces{v_1}$.
		This immediately satisfies~\labelcref{item:vanilla_version_verification_size,item:vanilla_version_verification_inclusion}.
		\labelcref{item:vanilla_version_verification_connected_in_time} holds true because every vertex in $\gm$ other than $v_1$ got $M$-connected to $v_1$ right after it was chosen as a fresh vertex.
		
		\textbf{Recursion step}:
		Now suppose $K^n$ and $W_n$ subject to~\labelcref{item:vanilla_version_verification_size,item:vanilla_version_verification_inclusion,item:vanilla_version_verification_connected_in_time} are given for some fixed $n \in \mbbN$.
		Consider the first $n+1$ vertices that are completely $M$-adjacent to every $v \in V(K^n)$.
		Such vertices exist since every $w$ in the infinite set $W_n$ has this property by~\labelcref{item:vanilla_version_verification_connected_in_time}.
		Let us call this set of vertices $F$.
		
		Now consider any vertex $w \in W_n \setminus F$.
		At the point in the game $n$ turns after $w$ was chosen as a fresh vertex by Maker, it was already $M$-connected to $V(K^n)$, and at that point $w$ was $B$-adjacent to at most $n$ other vertices, so at least one vertex $v^\prime \in F$ was still available.
		Let $\hat{v}$ be the vertex in $F$ with this property and the smallest possible index.
		Then Maker claimed $w\hat{v}$ in her $(n+1)^{\text{st}}$ move of $M$-connecting $w$.
		As $w$ was arbitrary, this is true for every one of the infinitely many vertices of $W_n \setminus F$ and so, as $F$ is finite, at least one vertex of $F$ gets chosen in this way for infinitely many vertices from $W_n$.
		We denote the smallest such vertex in $F$ by $v^\ast$ and set
		\begin{align*}
			K^{n+1} := & \left( V(K^n) \cup \braces{v^\ast}, E(K^n) \cup \braces{v v^\ast \colon v \in V(K^n)} \right), \text{ as well as}\\
			W_{n+1} := & \left\lbrace w \in W_n \colon v^\ast \text{ was } M\text{-connected in the }(n+1)\text{\textsuperscript{st} turn} \right.\\ 
			& \left. \text{after being picked as a fresh vertex} \right\rbrace.
		\end{align*}
		This takes care of~\labelcref{item:vanilla_version_verification_inclusion} and additionally, we have $\verts{K^{n+1}} = n+1$ and $\verts{W_{n+1}} = \aleph_0$, thus~\labelcref{item:vanilla_version_verification_size} is satisfied.
		\labelcref{item:vanilla_version_verification_connected_in_time} follows from the recursion assumption together with the choice of $W_{n+1}$.
	\end{proof}

	As we know now that the $\Kalephzero$-game is a Maker win, we will go on to consider some variants in which we make life a little harder for her.
	The natural way to do so might be to allow Breaker to claim more than one edge for any edge that Maker claims.
	This kind of variant has also been studied in the finite case, see~\cite{HKSS14}*{Chapter 3} and is called a \emph{biased} game.
	In our setting it doesn't make much difference, at least as long as Breaker is only allowed to claim the same finite number $k$ of edges on each turn.
	Maker does not even need to adapt her strategy.
	In the verification we need to take $F$ to be of size $k n+1$ rather than $n + 1$, and the argument works just as before.
	
	What if Breaker is allowed to claim a monotone increasing number of edges on his turns? It turns out that regardless of how slow the increment actually is, as long as the number of edges he claims tends to infinity he has a winning strategy:
	
	At the beginning of the game, he picks an enumeration $e_1 = \braces{x_1,y_1},e_2 = \braces{x_2,y_2},e_3= \braces{x_3,y_3},\dots$ of the edges of $E$.
	For any $n \in \mbbN$ there is an $N \in \mbbN$ such that  from the $N^\text{th}$ turn on, Breaker is allowed to claim $n$ edges in each of his turns, for any edge $e = \braces{x,y}$ that Maker claims.
	Beginning at $i=1$, for every $i \leq n$, whenever $G[\braces{x,y,x_i,y_i}] \not \subseteq \gm \cup \gb$, Breaker claims one of the available edges from $G [ \braces{x,y,x_i,y_i}] $ in his $i$\textsuperscript{th} turn.
	This strategy ensures that $e_n$ can be part of a complete graph in $\gm$ of at most some finite size dependent on $N$.
	As this holds for any $n \in \mbbN$, Maker cannot construct a $\Kalephzero \subseteq \gm$.
	
	In the following sections we will take a closer look at other, more challenging variations.
	One could also consider biased versions of the games considered later in this paper, but the theory of such biased games is always just the same as that outlined above and so we will not discuss it further.
	
	\section{Finitely many colours}\label{section:finitely_many_colours}
	A more interesting way to make Maker's objective more demanding is the following:
	Before the beginning of the game, every vertex of the board gets assigned one of $k \in \mbbN \setminus \braces{0}$ many colours.
	Let $k \in \mbbN$.
	A map 
	\[
	c \colon V(G) \longrightarrow [k]
	\]
	or
	\[
	c \colon V(G) \longrightarrow \mbbN
	\]
	is a \emph{colouring} of $V(G)$ if $\left\vert c^{-1}(i) \right\vert = \aleph_0$ for every colour $i \in [k]$.
	Moreover, for a set $W \subseteq V(G)$ we define $c[W] := \braces{c(v) \colon v \in W}$.
	Let $c$ be a colouring and $j \in \im (c)$.
	We call 
	\begin{displaymath}
		c^{-1}(j) \subseteq V(G)
	\end{displaymath}
	
	the \emph{colour class} of $j$.
	
	Maker's objective will be to build a $\Kalephzero \subseteq G_M$ as before but with the additional property that it includes infinitely many vertices from every colour class.
	We will call this version of the game the finitely coloured $\Kalephzero$-game.
	This is again a Makers' win with the following strategy.
	
	\medskip
	Let $k \in \mbbN \setminus \braces{0}$, suppose that the board is coloured by a colouring $c : V(G) \longrightarrow [k]$ and let $v_n$ be the vertex added to $\gm$ most recently.
	Now suppose $\deg_M(v_n) \equiv \ell \mod k$ and $v_n \in c^{-1}(h)$ for $h,\ell \in [k]$ not necessarily distinct.
	Then if Maker connects $v_n$ to a vertex of colour $\ell$ in the following fashion, we say that she plays according to the \emph{finite colour balanced greedy strategy}.
	
	Let $F \subseteq V(\gm)$ be the set of the first $k \cdot \deg_M(v_n) +1$ many vertices such that for all $v_m \in F$:
	\begin{itemize}
		\item $\nm (v_n) \subseteq \nm (v_m)$,
		\item $v_m \in c^{-1} (\ell)$,
		\item $m<n$, and
		\item $a  < m $ for all $v_a \in \nm (v_n)$.
	\end{itemize}
	If there are fewer than $k \cdot\deg_M(v_n) +1$ vertices satisfying these conditions, Maker chooses a fresh vertex $v_{n+1}$ of colour $m$, where $n+1 \equiv m \mod k$,  and claims $v_1v_{n+1}$.
	Otherwise, she considers the set $K \subseteq V(\gm)$ of all vertices $v_i $ satisfying:
	\begin{itemize}
		\item $j < i$ for all $v_j \in F$,
		\item $ \nm (v_n)  \supseteq \nm (v_i)$, and
		\item $v_i \in c^{-1}(h)$.
	\end{itemize}
	Maker assigns a tuple in $\mbbN \times \mbbN$ to every $v_i \in F$ via the injective map
	\begin{align*}
		f: F & \longrightarrow \mbbN \times \mbbN,\\
		v_i & \longmapsto \left( \verts{\nm (v_i) \cap K} , i \right)
	\end{align*}
	and then she orders $f(F)$ lexicographically, which results in an ordered set
	\begin{align*}\numberthis\label{equation:finitely_many_colourclasses_strategy_ordering}
		\left( f(F), \leq \right).
	\end{align*}
	Maker determines the smallest tuple $(\verts{\nm (v_{\delta}) \cap K},{\delta}) \in f (F)$ such that \linebreak ${v_{\delta}v_n \notin E(\gb)}$ and claims this edge.
	By the size of $F$ it is clear that there will be a vertex $v_{\delta}$ available, as Breaker had only $\deg_M(v_n) < k \cdot \deg_M(v_n) +1$ many moves where he could have coloured edges that are incident with $v_n$.
	\medskip
	
	Let us shed some light on two aspects of this strategy, namely the size of $F$ and the purpose of the order on $F$ induced by $f$.
	
	Our verification that this strategy works will be similar to that in \Cref{section_TheOriginalKOmegaGame}, in that we will again recursively build a nested sequence of complete graphs $K^n$ for every $n$ and in every step make sure that there is an infinite set $W_n \subseteq V(\gm)$ such that for every vertex $v \in W_n$ the entire $K^n$ is contained in its neighbourhood, i.e. the induced subgraph on $V(K^n) \cup \braces{v}$ is a potential candidate to continue the sequence.
	Then the crucial part is to carefully pick a vertex such that there still is an infinite set $W_{n+1} \subseteq W_n$ left.
	Note that the set $K$ for some vertex $v$ in the strategy will be contained in the corresponding set $W_n$ in the proof if $v$ is considered as a potential next vertex in the recursion.
	Because of the role the sets $W_n$ and therefore the sets $K$ play, we will informally refer to them  as \emph{reservoir}.
	In contrast to the proof in \Cref{section_TheOriginalKOmegaGame}, we need to also make sure that the sets $W_n$ also contain infinitely many vertices of every colour.
	This is precisely the motivation for the map $f$ introduced in the strategy above:
	If Maker just chose to play to the vertex from $ F$ with the smallest possible index as she does in the basic version, Breaker could ensure that all elements of the reservoir of colour $a$ are joined to some vertex $v_a$, but that all elements of the reservoir of colour $b$ are joined to some other vertex $v_b$.
	Then there would be no vertex that has infinitely many neighbours of both colours.
	Thus, instead of designating one vertex that has infinitely many neighbours of every colour, Maker instead ensures that for any colour $o$, Breaker can bar at most $\deg_M(v_n)$ vertices of $F$ from having infinitely many neighbours of colour $o$.
	This excludes at most $k \cdot m$ vertices of $F$ (recall that $k$ is the number of colours and $m$ the current $M$-degree of $v_n$).
	Maker wants to utilise this fact in order to ensure that there is at least one suitable vertex, i.e. a vertex with infinitely many neighbours of every colour, in the recursion step of the proof.
	She can achieve this by ensuring that the connection from vertices in the reservoir are spread as evenly as possible across $F$.
	The tool to do this is the function $f$ and the lexicographic ordering:
	
	Picking $v_{\delta}$ minimally in~\labelcref{equation:finitely_many_colourclasses_strategy_ordering} makes the choice of the vertex unique for Maker, this is ensured by the second entry of the ordering.
	More importantly, as we have argued above, the vertices of $F$ must be $M$-connected in a balanced fashion and this is achieved by choosing $v$ such that $\verts{\nm (v) \cap K}$ is smallest possible.
	To illustrate what we mean by that, one may think of the vertices in $K$ as being the set of vertices in $\gm$ that are identical to $v_n$ in the following sense:
	They were added to $\gm$ later than all of the vertices in $F$, the vertices got $M$-connected to $\gm$ during their first $\deg_M(v_n)$ many turns in the same manner as $v_n$, and they have the same colour as $v_n$, namely $h$.
	Via $f$, Maker finds the elements in $F$ that have the fewest neighbours in $K$ and out of these she chooses the one that has the smallest index.
	
	Therefore, by ensuring that $F$ has size $k \cdot m + 1$ playing to vertices of $F$ as evenly as possible via $f$ Maker ensures that there is a suitable vertex in the recursion step of the proof.
	
	\begin{theorem}\label{theo:finitely_many_colours}
		The finite colour balanced greedy strategy is a winning strategy for Maker in the finitely coloured version of the $\Kalephzero$-game.
	\end{theorem}
	\begin{proof}
		We want to show that at the end of the game, if Maker plays according to the strategy given, there is a $\Kalephzero \subseteq \gm$ that uses infinitely many vertices of each colour class.
		
		\textbf{Recursive construction:}
		For every $n \in \mbbN$ we will construct a complete graph $K^n \subseteq G_M$ together with a set of vertices $W_n \subseteq V(\gm)$ with the properties
		\begin{enumerate}[label={($\blacktriangle${\arabic{*}})}]
			\item\label{item:finitely_many_colourclasses_verification_inclusion} $K^{n}\subset K^{n+1}$,
			\item\label{item:finitely_many_colourclasses_verification_size_Wn} $\verts{W_n \cap c^{-1} (i)} = \aleph_0$ for all $i \in [k]$,
			\item\label{item:finitely_many_colourclasses_verification_connected_in_time} for each $w \in W_n$ we have $N_M(w) \supseteq V(K^n)$ and the vertices of $K^n$ were the first $n$ to become $M$-connected to $w$, and
			\item\label{item:finitely_many_colourclasses_verification_colours_of_Kn} $\verts{K^n} =n$ and there is an enumeration $\left\lbrace v^\prime_1, v^\prime_2, \dots , v^\prime_n \right\rbrace$ of $V(K^n)$ such that $v_j^\prime$ is coloured in $m$ and $j \equiv m \mod k$ for every $1\leq j \leq n$.
		\end{enumerate}
		Note that by~\labelcref{item:finitely_many_colourclasses_verification_size_Wn} we have in particular $\verts{W_n} = \aleph_0$.
		As in the proof of \Cref{theo:Original_Game}, we can get the desired $\Kalephzero$ from properties~\labelcref{item:finitely_many_colourclasses_verification_inclusion,item:finitely_many_colourclasses_verification_colours_of_Kn} by considering
		\[
		\bigcup_{n \in \mbbN} K^n = \Kalephzero.
		\]
		Here~\labelcref{item:finitely_many_colourclasses_verification_inclusion} ensures that there is the sequence $K^1 \subset K^2 \subset K^3 \subset \dots $ of complete graphs.
		\labelcref{item:finitely_many_colourclasses_verification_colours_of_Kn} ensures that there are infinitely many vertices of each colour.
		\labelcref{item:finitely_many_colourclasses_verification_size_Wn,item:finitely_many_colourclasses_verification_connected_in_time} are needed to ensure that there always is a next vertex that can be added to $K^n$ to form $K^{n+1}$.
		It remains to show that the conditions above can be preserved in every step.
		
		\textbf{Initial step:}
		Again, we can set $K^1 := \left( \lbrace v_1 \rbrace , \emptyset \right)$ and $W_1 = V(\gm) \setminus \braces{v_1} $.
		As this is the initial step,~\labelcref{item:finitely_many_colourclasses_verification_inclusion} holds true.
		Since Maker repeatedly added vertices of all colours to $\gm$,~\labelcref{item:finitely_many_colourclasses_verification_size_Wn} is true as well.
		\labelcref{item:finitely_many_colourclasses_verification_connected_in_time} holds, since $v_1$ was the first vertex to be joined to each $v_i$ with $i \in \mbbN \setminus \braces{1}$.
		Finally, as $V(K^1) = \braces{v_1}$ and $v_1$ is coloured with colour $1 \in [k]$,~\labelcref{item:finitely_many_colourclasses_verification_colours_of_Kn} is also true.
		This concludes the base case.
		
		\textbf{Recursion step:} 
		Let $n\in \mbbN$, $1\leq i \leq k$, let $K^n$ and $W_n$ subject to~\labelcref{item:finitely_many_colourclasses_verification_inclusion,item:finitely_many_colourclasses_verification_colours_of_Kn} be given and let $i \in [k]$ such that $n +1 \equiv i \mod k$.
		We want to construct $K^{n+1}$ and $W_{n+1}$ with the required properties.
		In order to do so, let $F$ be the set of the first $k n+1$ vertices of colour $i$ that have a common magenta edge with every vertex of $K^n$.
		Such a set exists, since all of the infinitely many vertices of colour $i$ in $W_n$ have this property.
		
		Let $j$ be the largest index of a vertex in $F$, fix an arbitrary colour $\ell \in [k]$ and let $w \in \left( W_n\setminus F \right) \cap c^{-1}(\ell)$ be a vertex with an index larger than $j$.
		After Maker claimed $wv^\prime$ for every $ v^\prime \in V(K^n)$ in her first $n$ moves of connecting $w$ to $\gm$, the statement
		\[
		\nm (w) \subseteq\nm (v)
		\]
		held for all $v \in F$.
		Since these are the first $kn+1$ such vertices, Maker chose the smallest available one of them with respect to the ordering derived from $f$ as defined in \labelcref{equation:finitely_many_colourclasses_strategy_ordering}.
		Breaker could block at most $n$ edges $vw$ for $v \in F$, thus there are at least $(k-1)n+1$ possible edges for Maker to choose from.
		Therefore, at most $k$ vertices of $F$ individually have only finitely many vertices in $W_n \cap c^{-1}(l)$ that chose them in the $(n+1)^{\text{st}}$ move of connecting to $\gm$.
		
		As the colour $l$ was arbitrary, the argument above holds true for each of the $k$ different colours.
		Therefore there is at least one vertex $u$ in $F$ that has infinitely many neighbours in every colour class in $W_n$.
		We choose the vertex $u^\ast \in F$ of these with the smallest index and let $K^{n+1}$ be the graph obtained from $K^n$ by adding $u^\ast$ and all edges from it to $K^n$.
		As $W_{n+1}$ we take the set of vertices in $W_n$ such that $u^\ast$ was the ($n+1$)\textsuperscript{st} vertex to which they were $M$-connected.
		This ensures~\labelcref{item:finitely_many_colourclasses_verification_inclusion}.
		Moreover, it means that $W_{n+1}$ contains infinitely many vertices of every colour class by the choice of $u^\ast$, therefore ensuring~\labelcref{item:finitely_many_colourclasses_verification_size_Wn}.
		The first $n+1$ vertices to be joined to any $w \in W_{n+1}$ were those of the $ K^{n+1}$ by the induction hypothesis and the construction of $u^\ast$, so we have~\labelcref{item:finitely_many_colourclasses_verification_connected_in_time}.
		Lastly, the fact that we considered only vertices of colour $ n+1 \equiv i \mod k$ for $F$, together with the assumption that~\labelcref{item:finitely_many_colourclasses_verification_colours_of_Kn} holds for $n$, ensures~\labelcref{item:finitely_many_colourclasses_verification_colours_of_Kn} for step $n+1$.
	\end{proof}
	
	\section{Infinitely many colours}\label{sec:infinitely_many_colourclasses}
	Next we consider what happens if there are infinitely many colours.
	As before, we will consider a surjective map $c \colon  V(G) \longrightarrow \mbbN$ and we will require $\left\lvert c^{-1} (i) \right\rvert = \aleph_0$ for every $i \in \mbbN$.
	
	If it is Maker's aim to construct a $\Kalephzero$ that uses every colour class infinitely often, she is doomed to fail, as there is a strategy for Breaker with which he can keep Maker from doing so.
	However, if it is Maker's aim to construct a $\Kalephzero$ that only uses infinitely many different colour classes, then there is a strategy with which she can secure this.
	We will first present Breaker's strategy for the first variant and after that we will give the strategy according to which Maker should play to win the second variant.
	
	\subsection{Using all colours of the board}\label{subsec:infinitely_many_vertex_colourclasses_Breaker_strategy}
	Our aim is to define a pairing strategy such that for every edge $e$ of the board there is some colour class $i$ of which Maker may use no vertices together with $e$.
	First note that there are countably infinitely many edges in a $\Kalephzero$ as $V(\Kalephzero)$ is countably infinite and the edges correspond to the two element subsets of $V$.
	
	\medskip
	Before the beginning of the game, Breaker picks an enumeration $e_1,e_2,e_3,\dots$ of all the edges of the board $G$.
	He then recursively finds an enumeration $c_1,c_2,c_3,\dots$ of infinitely many colours that are present in the $\Kalephzero$ such that for all $i \in \mbbN$:
	\begin{enumerate}[label=($\bigstar${\arabic{*}})]
		\item\label{no_previous_colour_recursion} $c_i \neq c_j$ for all $j<i$ and
		\item\label{no_colours_of_previous_vertices_recursion} $c_i \notin c \left[ \bigcup_{j \leq i} e_j \right]$.
	\end{enumerate}
	Fix some $m \in \mbbN$.
	We set $V_m := c^{-1}(c_m)$ and according to \labelcref{no_colours_of_previous_vertices_recursion} we have $V_m \cap e_m = \emptyset $.
	Suppose $e_m = xy$.
	Then for any $v \in c^{-1}(c_m) $ there are exactly two edges from $v$ to $e_m$, namely $vx$ and $vy$.
	Whenever Maker claims one of those edges, Breaker claims the other one in his following turn.
	\begin{lem}\label{lem:infinitely_many_colours_breaker_win}
		The all infinite colourclass pairing strategy is a pairing strategy and furthermore a winning strategy for Breaker in the infinitely coloured version of the $\Kalephzero$-game where Maker must have all colours contained in her $\Kalephzero$.
	\end{lem}
	\begin{proof}
		We first check whether the defined strategy actually is a pairing strategy, i.e. that any edge lies in at most one of the pairs of edges.
		Indeed, suppose for a contradiction, that there is an edge $e=uw$, that lies in two different pairs of edges.
		
		\textbf{Case 1}:
		$e$ lies in the pair of edges for two distinct edges $e_i$, $e_j$ that are incident with the same vertex of $e$, $u$ say:
		Then $w \in c^{-1} (c_i) \cap c^{-1} (c_j)$.
		Thus we have $c_i = c_j$, a contradiction to \labelcref{no_previous_colour_recursion}, as either $i<j$ or $j<i$.
		
		\textbf{Case 2}:
		$e$ lies in the pair of edges for two non adjacent edges $e_i$, $e_j$:
		Without loss of generality we may assume $u \in e_i$ and $w \in e_j$.
		This can only happen, if $u \in V_j$ and $w \in V_i$, a contradiction to \labelcref{no_colours_of_previous_vertices_recursion}, as either $i<j$ or $j<i$.
		
		Thus the pairs used for the strategy are disjoint.
		Therefore the given strategy actually is a pairing strategy for Breaker.
		
		Furthermore Maker cannot build a $\Kalephzero$ that uses all of the colours present on the board, if Breaker plays according to the defined strategy:
		For any edge that she wants to incorporate into her $\Kalephzero$, there is a colour class that corresponds to it according to the construction which she therefore cannot use in her $\Kalephzero$.
	\end{proof}
	Note that this result can be strengthened in the following sense:
	As every edge in a $\Kalephzero$ of Maker has a different colour assigned to it and the $\Kalephzero$ contains infinitely many edges, there are infinitely many colours of which Maker cannot incorporate infinitely many into her $\Kalephzero$, this means that Breaker can even stop Maker from using cofinitely many colours, each infinitely often, in a $\Kalephzero$.
	
	\subsection{Using infinitely many colours of the board}\label{subsec:infinitely_many_vertex_colourclasses_Maker}
	Let us now investigate how Maker should play in order to ensure that $\gm$ contains infinitely many vertices from infinitely many different colour classes.
	As before, she needs to add fresh vertices to $\gm$, making sure that she keeps track of the colours as well as taking care of the vertices that are already part of $\gm$, while also ensuring that they are as interconnected as possible in general and paying attention to the colours of these fresh vertices in particular.
	
	We first introduce one additional difinition.
	For a finite subset $U \subseteq V(\gm)$ we let
	\begin{align*}
		\varphi_U  \colon  \left[ \verts{U} \right] \longrightarrow \left\lbrace i \in \mbbN \colon  v_i \in U \right\rbrace\subseteq \mbbN
	\end{align*}
	be the unique order preserving bijection.
	For an infinite subset $W \subseteq V(\gm)$ we consider the unique order preserving bijection
	\begin{align*}
		\varphi_W  \colon  \mbbN \longrightarrow \left\lbrace i \in \mbbN \colon \  v_i \in W\right\rbrace\subseteq \mbbN.
	\end{align*}
	
	In this variant of the game, Maker cannot rotate through the colour classes in the same fashion as with finitely many colours and thus she has to work on her objective in a diagonal fashion.
	We further specify this in the strategy:
	
	At the beginning of the game, Maker chooses a sequence $s_1,s_2,s_3,\dots$ of all the colours appearing on the board such that each individual colour appears infinitely often.
	Let us call this sequence $S$.
	
	\medskip
	We call the following strategy for Maker the \emph{infinite colour balanced greedy strategy}.
	
	In her first turn, she picks two fresh vertices of colours $s_1$ and $s_2$, calls them $v_1$ and $v_2$ respectively and claims the edge $v_1v_2$ for herself.
	When Maker adds vertices to her subgraph in later stages of the game, letting $\verts{\gm} = n-1$, she adds a fresh vertex $v_{n}$ of colour $s_n$ to $\gm$.
	
	After $M$-connecting a fresh vertex $v_{n}$ of colour $s_n$ by claiming $v_1v_n$, on the next few turns Maker determines which edge $v_iv_n$ to claim by considering the set $U \subseteq V(\gm)$ of all vertices $v_u$ that satisfy
	\begin{itemize}
		\item $\nm (v_n) \subseteq \nm(v_u)$, moreover in the first $\deg_M(v_n)$ turns of $M$-connecting $v_u$ Maker $M$-connected $v_u$ to the same vertices as $v_n$, in the same order, and
		\item $i<u$ for all $v_i \in \nm (v_n)$.
	\end{itemize}
	Then Maker considers the subset 
	\begin{align*}\numberthis\label{equation:definition_of_U_prime}
		U^\prime := \left\lbrace v \in U \colon c(v) = c(v_n) \right\rbrace
	\end{align*}
	and with $\varphi_U$ as defined above she determines that the vertex she plays to next should be of colour
	\begin{align*}\numberthis\label{equation:definition_of_colour_j}
		j = c \left( v_{\varphi_{U} \left( \verts{U^\prime} \right)} \right).
	\end{align*}
	Note that $v_n \in U$, thus $\varphi_{U} \left( \verts{U^\prime} \right)$ is well defined.
	Next, Maker lets $F \subseteq V(\gm)$ be the set of the first $( \verts{c[\nm (v_n)]}+ 2) \cdot \deg_M(v_n) +1$ vertices, such that
	\begin{itemize}
		\item $k<i$ for all $v_i \in F$ and all $v_k \in \nm(v_n)$,
		\item $i<n$ for all $v_i \in F$,
		\item $v_i \in c^{-1} (j)$ for all $v_i \in F$, and
		\item all $v_i \in F$ satisfy $\nm (v_n) \subseteq \nm (v_i)$.
	\end{itemize}
	If there are fewer than $( \verts{c[\nm (v_n)]}+ 2) \cdot \deg_M(v_n) +1$ vertices satisfying these conditions, Maker instead chooses a fresh vertex $v_{n+1}$ as described above.
	If there is such a set, she chooses its subset of the first $( \verts{c[\nm (v_n)]}+ 2) \cdot \deg_M(v_n) +1$ many and calls this $F$.
	Maker wants to $M$-connect a vertex from $F$ to $v_n$ analogously to the finite colour balanced greedy strategy:
	she considers the set $K \subseteq V(\gm)$ of all vertices $v_k$ that satisfy
	\begin{itemize}
		\item for all $v_i \in F$ we have $i < k$,
		\item for all $v_\ell \in \nm (v_n)$ we have $\ell<k$,
		\item $\nm (v_k) \supseteq \nm (v_n)$, and
		\item $c(v_k) = c(v_n)$.
	\end{itemize}
	Maker assigns a tuple to every $v_i \in F$ as follows:
	\begin{align*}
		g \colon F & \longrightarrow \mbbN \times \mbbN,\\
		v_i & \longmapsto \left( \verts{\nm (v_i) \cap K} , i \right),
	\end{align*}
	and then orders $g(F)$ lexicographically, which results in an ordered set
	\begin{align*}\numberthis\label{equation:infinitely_many_colourclasses_strategy_ordering}
		\left( g(F), \leq \right).
	\end{align*}
	Maker determines the smallest $v_{\delta} \in F$ such that $v_{\delta}v_n \notin E(\gb)$ and claims this edge.
	\medskip
	
	Note that as $\verts{F} = ( \verts{c[\nm (v_n)]}+ 2) \cdot \deg_M(v_n) +1 $, there will be a vertex $v_{\delta}$ available, as Breaker has only had $\deg_M( v_n )$ many moves where he could have $B$-connected vertices from $F$ with $v_n$.
	As in \Cref{section:finitely_many_colours} considering the ordering $(g(F), \leq)$ ensures that Maker plays from vertices similar to $v_n$ to the vertices in $F$ in a ``balanced fashion'', which is crucial in the verification step.
	
	Let us investigate why the size of $F$ should be $( \verts{c[\nm (v_n)]}+ 2) \cdot \deg_M(v_n) +1$.
	Recall that the size of the corresponding set in the finite colour balanced greedy strategy was ``(number of colours on the board $\cdot$ degree of the active vertex)$ +1$''.
	The sizes thus only differ by ``number of colours on the board'' vs ``number of colours in the neighbourhood of the active vertex$ + 2$''.
	It is clear that ``number of colours on the board'' cannot be used in the infinitely coloured game, as the size of $F$ must be finite.
	It becomes clear, why the chosen number gives a good compromise for the following reason, when we suppose that the vertex $v_n$ will be considered as an element of the set $W_m$ of potential future vertices for some $m \in \mbbN$.
	In the proof we must make sure that for any colour $d$ already present on the $K^m$ there are infinitely many vertices of colour $d$ in $W_m$.
	This is ensured by the ``number of colours in the neighbourhood of the active vertex''-part.
	On top of that, in order to eventually have infinitely many colours present on the $\Kalephzero$, we need to (a) allow for one additional colour in case we want to add a new colour to the $K^{m+1}$ and (b) ensure that there will still be infinitely many other potential colours to add in the future present in $W_{m+1}$.
	This gives rise to the ``$+2$''-part.
	
	Lastly, let us elucidate Makers' choice of the colour $j$ given in \labelcref{equation:definition_of_colour_j}  before we move on.
	Breaker might render some colours unusable for Maker, but which these will be will not be clear until after the game.
	Thus, Maker needs a method to ensure that for each two such colours $j$ and $k$ she infinitely often tries to connect from a vertex of colour $k$ down to one of colour $j$. 
	The given function fulfils this purpose, which we will prove later (see \labelcref{equ:ProofThatColourpickingWorks}).
	
	\begin{theorem}\label{theo:infinitely_many_colours_maker_win}
		The infinite colour balanced greedy strategy is a winning strategy for Maker in the infinitely coloured version of the $\Kalephzero$-game where Maker must have infinitely many colours contained infinitely often in her $\Kalephzero$
	\end{theorem}
	\begin{proof}
		We want to prove that, after infinitely many turns, there is a $\Kalephzero \subseteq \gm$ that uses infinitely many vertices of infinitely many different colours, if Maker plays according to the strategy above.
		Before we begin with the recursion, we pick a sequence $\hat{C} = c_1,{c}_2,{c}_3,\dots$ of colours of $c[V]$ which contains every element of $c[V]$ infinitely often.
		We just require that ${c}_1 = s_1 = c(v_1)$.
		
		\textbf{Recursive construction:}
		For every $n \in \mbbN \setminus \braces{0}$ we will construct a complete graph $K^n \subseteq \gm$ together with a set of vertices $W_n \subseteq V$, and a set of colours $C_n \subseteq c[V]$ with the properties
		\begin{enumerate}[label={($\blacklozenge${\arabic{*}})}]
			\item\label{infinite_vertex_colouring_properties_of_K_n} $K^n \subset K^{n+1}$,
			\item\label{infinite_vertex_colouring_colours_present_in_W_n} $\verts{W_n \cap c^{-1 } (i)} = \aleph_0$ for every $i \in C_n$,
			\item\label{infinite_vertex_colouring_vertices_got_joined_fast} for each $w \in W_n$ the first $n$ moves of connecting $w$ to $\gm$ by Maker were claiming the edges that join $w$ to the $K^n$,
			\item\label{infinite_vertex_colouring_verification_colours_of_Kn} $\verts{K^n} = n$ and there is an enumeration $\braces{v_1^\prime, v_2^\prime, \dots , v_n^\prime}$ of $V(K^n)$ such that $c (v_i^\prime) = c(v_j^\prime)$ if and only if $c_i = c_j$ for $1 \leq i \leq j \leq n$, and
			\item\label{infinite_vertex_colouring_colour_reservoir} $\verts{C_n} = \aleph_0$ and $c [V(K^n)] \subseteq C_n$.
		\end{enumerate}
		Note that in \labelcref{infinite_vertex_colouring_verification_colours_of_Kn}~we do not require $c(v_i) = c_i$.
		This is indeed impossible to achieve.
		But it secures that any colour that appears really appears infinitely often and together with \labelcref{infinite_vertex_colouring_colour_reservoir} it furthermore secures that infinitely many different colours appear in the inclusive chain
		\[
		K^1 \subset K^2 \subset K^3 \subset K^4 \subset K^5 \subset \dots ~,
		\]
		thus
		\[
		\bigcup_{n \in \mbbN} K^n = \Kalephzero
		\]
		is the desired complete subgraph of $\gm$.
		
		\textbf{Initial step}:
		Set $K^1 := \left( \lbrace v_1 \rbrace , \emptyset \right)$, $W_1 := V(\gm) \setminus \braces{v_1}$ and $C_1 = c[W_1]$.
		\labelcref{infinite_vertex_colouring_properties_of_K_n} holds true, since this is the initial step.
		As $\verts{K_1} = 1$, \labelcref{infinite_vertex_colouring_verification_colours_of_Kn} holds true as well.
		$\verts{C_1} = \aleph_0$ is true and $c[V(K^1)] \subseteq C_1$ is satisfied because $S$ contains every colour infinitely often, thus $C_1$ satisfies~\labelcref{infinite_vertex_colouring_colour_reservoir}.
		Moreover, as every vertex of $\gm \setminus \braces{v_1}$ was first $M$-connected to $v_1$,~\labelcref{infinite_vertex_colouring_vertices_got_joined_fast} is true.
		Finally, as $S$ contains every colour infinitely often, this ensures~\labelcref{infinite_vertex_colouring_colours_present_in_W_n} for $C_1$.
		This concludes the base case.
		
		\textbf{Recursion step}:
		Let $n \geq 1$ and $K^n$, $W_n$ and $C_n$ subject to~\labelcref{infinite_vertex_colouring_properties_of_K_n,infinite_vertex_colouring_colour_reservoir} be given and suppose $c_m$ was the entry of $\hat{C}$ we worked with in the previous step.
		(This means in particular that $c(v_n) = c_m$.)
		We want to construct $K^{n+1}$, $W_{n+1}$ and $C_{n+1}$ with the required properties.
		As before, we will sometimes need to make sure that we add vertices of colours that are already present in the $K^n$ but we will sometimes also need to add vertices of colours that are not.
		If there is $i \leq n$ such that $c_i = c_{n+1}$, we set $c_p := c(v_i)$ and otherwise we choose $c_p \in C_n \setminus c[ V(K^n) ]$ arbitrarily.
		We want to add a vertex of colour $c_p$ next and let $F$ be the set of the first $\left( \verts{c[V(K^n)]} + 2 \right) \cdot n +1$ vertices of colour $\thiscolour$ that have a common magenta edge with every vertex of the $K^n$.
		This set exists since there are infinitely many vertices of colour $\thiscolour$ in $W_n$ by~\labelcref{infinite_vertex_colouring_colours_present_in_W_n}.
		Moreover, we need to restrict $W_n$ to only contain vertices whose $(n+1)^{\text{st}}$ turn of connecting it to $\gm$ was a vertex of colour $\thiscolour$ and we want to ensure that~\labelcref{infinite_vertex_colouring_colours_present_in_W_n} holds for this restriction of $W_n$ as well.
		Note that since~\labelcref{infinite_vertex_colouring_properties_of_K_n,infinite_vertex_colouring_verification_colours_of_Kn,infinite_vertex_colouring_colour_reservoir} are independent of $W_n$, they still hold and~\labelcref{infinite_vertex_colouring_vertices_got_joined_fast} will hold for the restriction as it is a subset of $W_n$.
		
		Fix an order preserving map $\psi \colon  \mbbN \setminus \braces{0} \to I$ such that
		\begin{align}\label{equ:ProofThatColourpickingWorks}
			\left\lbrace v_{\psi (i)} \colon i \in \mbbN \setminus \braces{0} \right\rbrace = W_n \cap c^{-1} (\thiscolour).
		\end{align}
		
		Then, for every $m \in \mbbN \setminus \braces{0} $ and every $d \in C_n$ the $\left( \psi (m) \right)$\textsuperscript{th} vertex of colour $d$ in $W_n$ got $M$-connected to a vertex of colour $\thiscolour$ in the $(n+1)$\textsuperscript{st} move of connecting it to $\gm$.
		Thus, there are infinitely many vertices of colour $d$ in $W_n$ whose $(n+1)$\textsuperscript{st} neighbour in $\gm$ (according to the order in which they were connected to it) was a vertex of colour $\thiscolour$.
		As $d$ was arbitrary, this is true for every colour in $C_n$.
		We can thus restrict $W_n$ to these vertices and work with this set $W^\prime_n$ from here on.
		
		Let $\ell \in c[V(K^n)] \cup \braces{\thiscolour}$.
		Since Maker played to the vertices of $F$ in a balanced fashion, there are at most $n$ vertices $v \in F$ such that only finitely many vertices $w \in W_n^\prime \cap c^{-1} (\ell)$ got $M$-connected to $v$ in their $(n+1)$\textsuperscript{st} move of connecting them to $\gm$.
		As $\ell$ was arbitrary, this is true for every colour in $ c[V(K^n)] \cup \braces{\thiscolour}$ and thus there are at least $n+1$ vertices that have infinitely many such vertices in $W_n^\prime \cap c^{-1} (p)$ for every $p \in c[V(K^n)] \cup \braces{\thiscolour}$.
		Conversely, regarding the infinitely many colours in $C_n \setminus ( c[V(K^n)] \cup \thiscolour)$, since Breaker can block at most $n$ vertices for any of them, there are at most $n$ vertices in $F$ that are chosen by only finitely many vertices of cofinitely many colours not yet occurring in the $K^n$ in the $(n+1)$\textsuperscript{st} move of connecting them to $\gm$.
		Combining this means that there is at least one vertex $u^\prime \in F$ that got chosen by infinitely many vertices of every colour in $c[V(K^n)] \cup \braces{\thiscolour}$ as well as infinitely many vertices of infinitely many distinct colours in $C_n$ in their $(n+1)$\textsuperscript{st} move of connecting them to $\gm$.
		We choose the smallest such vertex and call it $v^\prime_{n+1}$.
		We set
		\begin{itemize}
			\item $K^{n+1} := G \left[ V(K^n) \cup \lbrace v^\prime_{n+1} \rbrace \right]$,
			\item $C_{n+1}$ the set of colours $i \in C_n$ for which infinitely many vertices that lie in $c^{-1} (i) \cap W_n$ got $M$-connected to $v^\prime_{n+1}$ in their $(n+1)$\textsuperscript{st} move of connecting them to $\gm$, and
			\item $W_{n+1} \subseteq W_n^\prime$ as the vertices in $W_n$ that got $M$-connected to $v^\prime_{n+1}$ in their $(n+1)$\textsuperscript{st} move of connecting them to $\gm$ and that are coloured with a colour in $C_{n+1}$. 
		\end{itemize}
		This ensures~\labelcref{infinite_vertex_colouring_properties_of_K_n,infinite_vertex_colouring_colours_present_in_W_n,infinite_vertex_colouring_colour_reservoir}.
		\labelcref{infinite_vertex_colouring_verification_colours_of_Kn} holds true by the choice of $c_p$ and the definition of $v^\prime_{n+1}$.
		All vertices of $W_{n+1}$ are completely ($M$-)adjacent to the $K^n$ by the induction hypothesis and to $v^\prime_{n+1}$ according to the construction, so all vertices of $W_{n+1}$ are completely ($M$-)adjacent to the $K^{n+1}$.
		It follows from the induction hypothesis and the choice of $v^\prime_{n+1}$ that those were the first $n+1$ moves that Maker made for each element of $W_{n+1}$.
		This verifies~\labelcref{infinite_vertex_colouring_vertices_got_joined_fast}.
		
		This shows that all of the required properties are preserved throughout the induction and thus the claim is proved.
	\end{proof}

	\section{Open problems}\label{sec:open_problems}
	One obvious variation of the game proposed in \Cref{section_TheOriginalKOmegaGame} immediately comes to mind: while in \Cref{section:finitely_many_colours,sec:infinitely_many_colourclasses} the basic version of the game was altered by colouring the vertices, one could instead colour the edges of the board and again demand that $\gm$ contains an isomorphic copy of the board as a subgraph.
	We will call this game the $\Kalephzero$ \emph{edge colouring game}.
	While in the vertex case the colour classes are very symmetric, this is different for edge colourings: it could happen that the subgraph induced by some colour class is locally finite, while the one for another is not.
	\begin{question}
		Let $n \in \mbbN$. For which colourings $c \colon E(\Kalephzero) \longrightarrow [n] $ is there a winning strategy for one of the players in the $\Kalephzero$ edge colouring game?
	\end{question}
	Certainly one may consider a colouring with infinitely many colours as well.
	\begin{question}
		For which colourings $c \colon E(\Kalephzero) \longrightarrow \mbbN $ is there a winning strategy for one of the players in the $\Kalephzero$ edge colouring game?
	\end{question}
	Another possible generalisation is to adapt the game to hypergraphs.
	\begin{question}
		Let $H$ be a complete infinite $k$-regular hypergraph.
		In the Maker-Breaker game on $H$ where it is Makers aim to have an isomorphic copy of $H$ be contained in $\gm$, is there a winning strategy for Maker?
	\end{question}
	Naturally, we can consider a $\leq$$k$-regular hypergraph or even an infinite complete hypergraph instead of a $k$-regular one.
	Furthermore, we could also apply a vertex or an edge colouring to the board in any of these variants, just as in the $\Kalephzero$-game.
	Since the $k$-regular hypergraph variant is already very advanced, we will not state these as questions here and rather highlight one particularly intriguing question.
	\begin{question}
		Consider the Maker-Breaker game in which the players alternately claim finite subsets of $\mbbN$ and it is Makers aim to claim an infinite set $\mathcal{F}$ of pairwise disjoint finite subsets of $\mbbN$ as well as the union of every finite subset of $\mathcal{F}$.
		Is there a winning strategy for one of the players?
	\end{question}
	
	\appendix
	\section{Determinacy of Maker-Breaker games}\label{sec:Appendix}
	Although Maker-Breaker games are easier to handle than infinite games in general, we show here that they are still not necessarily determined. More precisely, assuming the continuum hypothesis we construct an infinite Maker-Breaker game which is not determined. 
	
	Let $X$ be any countably infinite set. Let $(\sigma_i : i < \omega_1)$ be an enumeration of the possible strategies for Maker on this set and $(\tau_i : i < \omega_1)$ an enumeration of the possible strategies for Breaker. We will construct sequences $(S_i : i < \aleph_1)$ and $(T_i : i < \aleph_1)$ of infinite subsets of $X$ recursively such that each $S_i$ meets each $T_j$, and so that Maker may claim $S_i$ if Breaker plays according to $\tau_i$ and Breaker may claim $T_i$ if Maker plays according to $\sigma_i$. 
	
	Suppose that all $S_j$ and $T_j$ with $j < i$ have already been constructed. Let $(U_n : n < \omega)$ be a sequence of infinite subsets of $X$ in which all $T_j$ with $j < i$ appear. Then we take $S_i$ to be the set of elements claimed by Maker in a play of the game in which Breaker plays according to $\tau_i$ and Maker, on her $n$\textsuperscript{th} move, claims some as yet unclaimed element of $U_n$. The construction of $T_i$ is similar.
	
	Now consider the Maker-Breaker game on $X$ whose set of winning sets is $\{S_i : i < \omega_1\}$. No $\tau_i$ can be a winning strategy for Breaker, since Maker can still claim $S_i$ and win if Breaker plays according to $\tau_i$. Similarly, no $\sigma_i$ can be a winning strategy for Maker, since Breaker can still claim $T_i$, which meets all $S_j$, and so prevent Maker from winning if Maker plays according to $\sigma_i$. So this game is not determined.
	
	\medskip
	\bibliography{ref}

\begin{bibdiv}
\begin{biblist}

\bib{B02}{article}{
      author={Beck, József},
       title={Ramsey games},
        date={2002},
        ISSN={0012-365X},
     journal={Discrete Mathematics},
      volume={249},
      number={1},
       pages={3\ndash 30},
  url={https://www.sciencedirect.com/science/article/pii/S0012365X01002242},
        note={Combinatorics, Graph Theory, and Computing},
}

\bib{B08}{book}{
      author={Beck, J{\'o}zsef},
       title={Combinatorial games: tic-tac-toe theory},
   publisher={Cambridge University Press Cambridge},
        date={2008},
      volume={114},
}

\bib{D17}{book}{
      author={Diestel, Reinhard},
       title={Graph theory},
     edition={5},
   publisher={Springer},
        date={2018},
        ISBN={978-3-662-53621-6},
}

\bib{E19}{misc}{
      author={Erde, Joshua},
       title={personal communication},
        date={2019},
}

\bib{HKNPRS17}{article}{
      author={Hefetz, Dan},
      author={Kusch, Christopher},
      author={Narins, Lothar},
      author={Pokrovskiy, Alexey},
      author={Requilé, Clément},
      author={Sarid, Amir},
       title={Strong ramsey games: Drawing on an infinite board},
        date={2017},
        ISSN={0097-3165},
     journal={Journal of Combinatorial Theory, Series A},
      volume={150},
       pages={248\ndash 266},
  url={https://www.sciencedirect.com/science/article/pii/S0097316517300353},
}

\bib{HKSS14}{book}{
      author={Hefetz, Dan},
      author={Krivelevich, Michael},
      author={Stojaković, Miloš},
      author={Szabó, Tibor},
       title={Positional {Games}},
    language={en},
   publisher={Springer Basel},
     address={Basel},
        date={2014},
        ISBN={978-3-0348-0824-8 978-3-0348-0825-5},
         url={http://link.springer.com/10.1007/978-3-0348-0825-5},
}

\bib{L08}{misc}{
      author={Leader, Imre~B},
       title={Lecture on hypergraph games},
        date={2008},
}

\bib{M09}{book}{
      author={Moschovakis, Yiannis~N},
       title={Descriptive set theory},
     edition={2},
   publisher={American Mathematical Society},
        date={2009},
      number={155},
        ISBN={978-0-8218-4813-5},
}

\end{biblist}
\end{bibdiv}
\end{document}